\documentclass[a4paper,12pt]{amsart}
\usepackage{amssymb,amsmath}
\usepackage{amsfonts}
\usepackage{color,cite}
\usepackage[normalem]{ulem}

\newcommand{\R}{\Bbb{R}}

\newcommand{\n}{\noindent}

\newcommand{\dis}{\displaystyle}

\newcommand {\eps}{\varepsilon}

\newtheorem{teorema}{Theorem}
\theoremstyle{plain}

\newtheorem{defi}{Definition}

\numberwithin{equation}{section}

\begin{document}
\title{A note on borderline Brezis-Nirenberg type problems}
\author{Julian Haddad}
\address{Departamento de Matem\'{a}tica \\
Universidade Federal de Minas Gerais\\
Belo Horizonte, 30123-970 Brazil}
\email{julianhaddad@gmail.com}
\urladdr{}
\thanks{}
\author{Marcos Montenegro}
\address{Departamento de Matem\'{a}tica \\
Universidade Federal de Minas Gerais\\
Belo Horizonte, 30123-970 Brazil}
\email{montene@mat.ufmg.br}
\urladdr{}
\thanks{The second author was partially supported by Capes, CNPq and Fapemig.}
\date{\today }

\begin{abstract}
This note concerns the existence of positive solutions for the boundary problems

\[
\left\{ \begin{array}{rcll}\dis
    - \mathcal{L}_A u & = & u^{2^* - 1} + \lambda u & \mbox{in} \ \Omega
    \\ u & = & 0 & \mbox{on} \ \partial\Omega\\
     \end{array}\right.,
\]

\[
\left\{ \begin{array}{rcll}\dis
    - \mathcal{L}_{a,p} u & = & u^{p^* - 1} + \lambda u^{p-1} & \mbox{in} \ \Omega
    \\ u & = & 0 & \mbox{on} \ \partial\Omega\\
     \end{array}\right.,
\]

\n where $\mathcal{L}_A u = {\rm div}(A(x) \nabla u)$ and $\mathcal{L}_{a,p} u = {\rm div}(a(x) |\nabla u|^{p-2} \nabla u)$ are, respectively, linear and quasilinear uniformly elliptic operators in divergence form in a non-smooth bounded open subset $\Omega$ of $\R^n$, $1 < p < n$, $p^* = np/(n-p)$ is the critical Sobolev exponent and $\lambda$ is a real parameter. Both problems have been quite studied when the ellipticity of $\mathcal{L}_A$ and $\mathcal{L}_{a,p}$ concentrate in the interior of $\Omega$. We here focus on the borderline case, namely we assume that the determinant of $A(x)$ has a global minimum point $x_0$ on the boundary of $\Omega$ such that $A(x) - A(x_0)$ is locally comparable to $|x - x_0|^\gamma$ in the bilinear forms sense to $|x - x_0|^\gamma I_n$, where $I_n$ denotes the identity matrix of order $n$. Similarly, we assume that $a(x)$ has a global minimum point $x_0$ on the boundary of $\Omega$ such that $a(x) - a(x_0)$ is locally comparable to $|x - x_0|^\sigma$. We provide a linking between the exponents $\gamma$ and $\sigma$ and the order of singularity of the boundary of $\Omega$ at $x_0$ so that these problems admit at least one positive solution for any $\lambda \in (0, \lambda_1(-\mathcal{L}_A))$ and $\lambda \in (0, \lambda_1(-\mathcal{L}_{a,p}))$, respectively, where $\lambda_1$ denotes the first Dirichlet eigenvalue of the corresponding operator.
\end{abstract}

\maketitle

\section{Introduction and statements}

A lot attention has been paid to a number of counterparts of the well-known Brezis-Nirenberg problem studied in \cite{BN}, which consists in determining the range of values $\lambda$ such that the problem

\begin{equation} \label{BN}
\left\{
\begin{array}{rcll}
- \Delta u &=& u^{2^* - 1} + \lambda u & {\rm in} \ \ \Omega,\\
u &=& 0  & {\rm on} \ \ \partial\Omega
\end{array}
\right.
\end{equation}

\n admits at least one positive solution, where $\Omega$ is a bounded open subset of $\R^n$, $n \geq 3$ and $2^* = 2n/(n-2)$ is the critical Sobolev exponent associated to critical embedding of the space $H^{1}_0(\Omega)$.

According to \cite{BN}, the problem (\ref{BN}) admits at least one positive solution for any $\lambda \in (0, \lambda_1(-\Delta))$ provided that $n \geq 4$, where $\lambda_1(-\Delta)$ denotes the first eigenvalue of the Laplace operator on $\Omega$ under Dirichlet boundary condition. Moreover, the problem admits no such a solution for $n \geq 3$ if either $\lambda \geq \lambda_1(-\Delta)$ or $\lambda \leq 0$ and $\Omega$ is a star-shaped $C^1$ domain. Regarding that $n = 3$ and $\Omega$ is a ball, a positive solution of (\ref{BN}) exists if, and only if, $\lambda \in (\frac{1}{4}\lambda_1(-\Delta), \lambda_1(-\Delta))$.

A closely related problem is

\begin{equation} \label{BN-p}
\left\{
\begin{array}{rcll}
- \Delta_p u &=& u^{p^*-1} + \lambda |u|^{p-2}u & {\rm in} \ \ \Omega\; ,\\
u &=& 0  & {\rm on} \ \ \partial\Omega
\end{array}
\right.
\end{equation}

\n where $\Delta_p u = {\rm div} ( |\nabla u|^{p-2} \nabla u )$ and $p^* = np/(n-p)$ denotes respectively the $p$-Laplace operator and the critical Sobolev exponent associated to critical embedding of the space $W^{1,p}_0(\Omega)$ for $1 < p < n$.

Some existence and nonexistence results are well known. Part of the results for $p = 2$ was extended by Egnell \cite{HE1}, Garcia Azorero and Peral Alonso \cite{AzPe} and Guedda and Veron \cite{GuVe}, who proved that the problem (\ref{BN-p}) admits at least one positive solution for any $\lambda \in (0, \lambda_1(-\Delta_p))$ provided that $n \geq p^2$, and admits no such a solution for $n > p$ if either $\lambda \geq \lambda_1(-\Delta_p)$ or $\lambda \leq 0$ and $\Omega$ is a star-shaped $C^1$ domain, where $\lambda_1(-\Delta_p)$ stands for the first eigenvalue of the $p$-Laplace operator on $\Omega$ under Dirichlet boundary condition. The case $p < n < p^2$ was studied by Egnell \cite{HE1}, who found a number $\lambda^* > 0$ such that a positive solution of (\ref{BN-p}) exists for any $\lambda \in (\lambda^*, \lambda_1(-\Delta_p))$. We also mention some existence results of nontrivial $C^1$ solutions which have been established when $\lambda \geq \lambda_1(-\Delta_p)$ in \cite{CFP} and \cite{GaRu} for $p=2$ and in \cite{ArGa} and \cite{DeLa} for $p \neq 2$.

We consider the following variants of the Brezis-Nirenberg problem

\begin{equation} \label{BN-g}
\left\{ \begin{array}{rcll}\dis
    - \mathcal{L}_A u & = & u^{2^* - 1} + \lambda u & \mbox{in} \ \Omega,
    \\ u & = & 0 & \mbox{on} \ \partial\Omega\\
     \end{array}\right.,
\end{equation}

\begin{equation} \label{BN-p-g}
\left\{ \begin{array}{rcll}\dis
    - \mathcal{L}_{a,p} u & = & u^{p^* - 1} + \lambda u^{p-1} & \mbox{in} \ \Omega,
    \\ u & = & 0 & \mbox{on} \ \partial\Omega\\
     \end{array}\right.,
\end{equation}

\n where $\Omega$ is a bounded open subset of $\R^n$ and $\mathcal{L}_A u = {\rm div}(A(x) \nabla u)$ and $\mathcal{L}_{a,p} u =$ ${\rm div}(a(x) |\nabla u|^{p-2} \nabla u)$ denote, respectively, linear and quasilinear operators of second order.

We assume that $A(x)$ is a positive definite symmetric matrix for each $x \in \overline{\Omega}$ with continuous entries on $\overline{\Omega}$ and $a(x)$ is a positive continuous function on $\overline{\Omega}$. So, both operators $\mathcal{L}_A$ and $\mathcal{L}_{a,p}$ are uniformly elliptic in $\Omega$ under the divergence form.

We are concerned here with existence of positive solutions for (\ref{BN-g}) and (\ref{BN-p-g}), namely solutions in the weak sense, respectively, in $H^1_0(\Omega)$ and $W^{1,p}_0(\Omega)$ and of $C^1$ class and positive in $\Omega$. Such solutions will be referred simply as positive $C^1$ solutions.

Denote by $\lambda_1(-\mathcal{L}_A)$ and $\lambda_1(-\mathcal{L}_{a,p})$ the first eigenvalues associated, respectively, to the operators $\mathcal{L}_A$ and $\mathcal{L}_{a,p}$ on $\Omega$ under Dirichlet boundary condition. By using an usual argument based on maximum principles, one easily follows that (\ref{BN-g}) and (\ref{BN-p-g}) admit no positive $C^1$ solution for $\lambda \geq \lambda_1(-\mathcal{L}_A)$ and $\lambda \geq \lambda_1(-\mathcal{L}_{a,p})$, respectively.

The problem of existence of positive $C^1$ solutions becomes interesting when the parameter $\lambda$ is below the associated first eigenvalue, in which operators and domains play an important role. Remark that the functions $\det(A(x))$ and $a(x)$ assume their global minimum values on $\overline{\Omega}$. In the case when some minimum point is in the interior of $\Omega$, existence and nonexistence of positive $C^1$ solutions of (\ref{BN-g}) and (\ref{BN-p-g}) have been addressed, respectively, in \cite{HE, MoMo2} and \cite{HE1, MoMo1}. Particularly, it has been proved in \cite{HE} and \cite{MoMo2} that the problem (\ref{BN-g}) admits at least one positive $C^1$ solution for any $\lambda \in (0, \lambda_1(-\mathcal{L}_A))$ provided that there exist constants $\gamma > 2$ and $C_0 > 0$ such that

\begin{equation} \label{H1}
A(x) \leq A(x_0) + C_0 |x - x_0|^\gamma I_n
\end{equation}
locally around a global minimum point $x_0 \in \Omega$ of $\det(A(x))$ in the sense of bilinear forms, where $I_n$ denotes the identity matrix of order $n$. Note that the condition (\ref{H1}) implies that the global minimum point $x_0$ of the determinant of $A(x)$ has order $\gamma$. The corresponding result for the problem (\ref{BN-p-g}) (see \cite{MoMo1}) states that at least one positive $C^1$ solution exists for any $\lambda \in (0, \lambda_1(-\mathcal{L}_{a,p}))$ by assuming the existence of a global minimum point $x_0 \in \Omega$ of $a(x)$ and constants $\sigma > p$ and $C_0 > 0$ such that

\begin{equation} \label{H2}
a(x) \leq a(x_0) + C_0 |x - x_0|^\sigma
\end{equation}
for all $x$ in a neighborhood of $x_0$. Both conditions (\ref{H1}) and (\ref{H2}) are sharp in the sense that the existence ranges decreases for lower orders $\gamma$ and $\sigma$, see \cite{MoMo2} and \cite{MoMo1}.

The case when $x_0$ is on the boundary of $\Omega$ has not been discussed so far. It is natural to hope that the shape (or regularity) of the border around a minimum point should exert some influence on this situation. Moreover, according to nonexistence results in the interior case, which easily extend to the borderline context (see for example Proposition 6 of \cite{HE} and Proposition 7 of \cite{HE1}), existence of positive $C^1$ solutions for any $\lambda \in (0, \lambda_1(-\mathcal{L}_A))$ or $\lambda \in (0, \lambda_1(-\mathcal{L}_{a,p}))$ leads naturally to lower bounds on the values of $\gamma$ and $\sigma$ as above.

In order to state our existence theorems, we introduce the following notion of singularity:

\begin{defi}
The boundary of an open subset $\Omega$ of $\R^n$ is said to be interior $\alpha$-singular at the point $x_0$, with $\alpha \geq 1$, if there exist a constant $\delta > 0$ and a sequence $(x_j) \subset \Omega$ such that $x_j \to x_0$ as $j \rightarrow + \infty$ and $B(x_j, \delta |x_j - x_0|^\alpha) \subseteq \Omega$.
\end{defi}

\n It is clear that the boundaries of $C^1$ domains or domains satisfying the interior cone condition are interior $1$-singular at each one of their points.

Statements below present connections in the borderline case between the exponents $\gamma$ and $\sigma$, given in (\ref{H1}) and (\ref{H2}), and the kind of singularity of the boundary of $\Omega$ at the point $x_0$, targeting to achieve existence for a boarder range on $\lambda$.

\begin{teorema}\label{T.1} Let $\Omega$ be a bounded open subset of $\R^n$ with $n \geq 5$. Assume that the determinant of $A(x)$ attains its global minimum at a point
$x_0$ on the boundary of $\Omega$ such that the comparison (\ref{H1}) is satisfied with $\gamma > \frac{2n - 4}{n - 4}$. Assume also that the boundary of $\Omega$ is interior $\alpha$-singular at $x_0$, with $\alpha \in [1, \gamma \frac{n - 4}{2n - 4})$. Then, the problem (\ref{BN-g}) admits at least one positive $C^1$ solution for any $\lambda \in (0, \lambda_1(-\mathcal{L}_A))$.
\end{teorema}

\begin{teorema}\label{T.2} Let $\Omega$ be a bounded open subset of $\R^n$ with $n> p^2$, where $p > 1$. Assume that the function $a(x)$ attains its global minimum at a point
$x_0$ on the boundary of $\Omega$ such that the comparison (\ref{H2}) is satisfied with $\sigma > \frac{np - p^2}{n - p^2}$. Assume also that the boundary of $\Omega$ is interior $\alpha$-singular at $x_0$, with $\alpha \in [1, \sigma \frac {n-p^2}{np - p^2})$. Then, the problem (\ref{BN-p-g}) admits at least one positive $C^1$ solution for any $\lambda \in (0, \lambda_1(-\mathcal{L}_{a,p}))$.
\end{teorema}

Remark that the assumptions $\gamma > \frac{2n - 4}{n - 4}$ and $\sigma > \frac{np - p^2}{n - p^2}$ are always sufficient conditions for $C^1$ domains or domains satisfying the interior cone condition and also imply that $\gamma > 2$ and $\sigma > p$. In addition, both theorems say that as greater are the exponents $\gamma$ and $\sigma$, respectively in (\ref{H1}) and (\ref{H2}), more interior singular is permitted the boundary of $\Omega$ to be at $x_0$.

Problems (\ref{BN}) and (\ref{BN-p}) should be faced with the interior context, once the matrix $A(x)$ and the function $a(x)$ are constant and results hold for general domains. Theorems \ref{T.1} and \ref{T.2} focus strictly on ellipticity concentrating only on the boundary of $\Omega$. Two direct prototypes of operators are provided by $A(x) = A_0 + C_0 |x - x_0|^\gamma I_n$ and $a(x) = a_0 + C_0 |x - x_0|^\sigma$, where $x_0 \in \partial\Omega$, $A_0$ is a positive semi-definite symmetric matrix of order $n$ and $\gamma$, $\sigma$, $a_0$ and $C_0$ are positive constants.

As usual, the proof of Theorems \ref{T.1} and \ref{T.2} is based on constrained minimization arguments. Minimax theorems can alternatively be employed in order. In any case, the key point consists generally in choosing a suitable family of test functions so that the asymptotic behavior of associated energy levels is comparable to best Sobolev constants. Unfortunately, the test functions used in the classical context and also in the interior case do not work when the basis point $x_0$ is on the boundary of $\Omega$. We recall below the proof outline via minimization and comment on a simple and efficient idea of construction of test functions (bubbles) to be used in the proof of the above results.

Precisely, one considers the functionals

\[
\Phi_{A}(u) = \int_\Omega \nabla^t u A(x) \nabla u\; dx - \lambda \int_\Omega u^2\; dx
\]

\n constrained to the set $E : = \{u \in W^{1,2}_0(\Omega):\; \int_\Omega |u|^{2^*}\; dx = 1\}$, where
$\nabla^t u$ denotes the transposed vector of $\nabla u$, and

\[
\Psi_{a,p}(u) = \int_\Omega a(x) |\nabla u|^p\; dx - \lambda \int_\Omega |u|^p\; dx
\]

\n constrained to the set $F : = \{u \in W^{1,p}_0(\Omega):\; \int_\Omega |u|^{p^*}\; dx = 1\}$. Note that $\Phi_{A}(u) = \Phi_{A}(|u|)$ for all $u \in E$ and $\Psi_{a,p}(u) = \Psi_{a,p}(|u|)$ for all $u \in F$, so by using the elliptic PDEs theory one easily checks that minimizers of $\Phi_{A}|_E$ and $\Psi_{a,p}|_F$ are positive $C^1$ solutions of (\ref{BN-g}) and (\ref{BN-p-g}), respectively.

According to Lemmas 2.1 of \cite{MoMo2} and \cite{MoMo1} (which also hold for $x_0$ in $\overline{\Omega}$), the strict inequalities

\begin{equation} \label{M1}
c_A := \inf_{u\in E} \Phi_{A}(u)  < m_A^{1/n} K(n,2)^{-2}
\end{equation}
and

\begin{equation} \label{M2}
c_{a,p} := \inf_{u\in F} \Psi_{a,p}(u)  < m_{a,p} K(n,p)^{-p}
\end{equation}
imply the existence of minimizers of $\Phi_{A}$ and $\Psi_{a,p}$, respectively, in $E$ and $F$, where $m_A$ and $m_{a,p}$ denote the global minimum values of $\det{A(x)}$ and $a(x)$ on $\overline{\Omega}$, respectively, and $K(n,p)^p$ denotes the sharp Sobolev constant for the classical Sobolev inequality on $C^1_0(\R^n)$, namely

\[
\left(\int_{\R^n} |u|^{p^*}\; dx\right)^{p/p^*} \leq K(n,p)^p \int_{\R^n} |\nabla u|^p\; dx\, .
\]
Note that these minimizes only are positive $C^1$ solutions, module a scaling, of the problems (\ref{BN-g}) and (\ref{BN-p-g}), provided that $c_A$ and $c_{a,p}$ are positive, respectively. On the other hand, the positivity of $c_A$ and $c_{a,p}$ is equivalent to the assumptions $\lambda < \lambda_1(-\mathcal{L}_A)$ and $\lambda < \lambda_1(-\mathcal{L}_{a,p})$ assumed in each one of theorems. Therefore, the conclusion of both theorems follows if we are able to prove strict inequalities (\ref{M1}) and (\ref{M2}) for any $\lambda > 0$. The proof of these inequalities bases on the construction of bubbles involving extremal functions for $K(n,p)^p$ in a similar way to the well-known works \cite{BN} and \cite{AzPe}, except that the mass of these functions should be concentrated in balls centered at interior points of $\Omega$ converging fast to $x_0$. The appropriate chosen of the speed of convergence of points and radius is the main part of proof. In the second section, we provide the needed bubbles estimates. In the third section, we prove (\ref{M1}) and (\ref{M2}) by using such estimates along with the assumptions assumed in each one of theorems.

\section{Bubbles estimates}

Let $\Omega$ be a bounded open subset of $\R^n$ and $x_0 \in \partial \Omega$. Assume that the boundary of $\Omega$ is $\alpha$-singular, with $\alpha \geq 1$, at $x_0$ in the sense of Definition 1.1. So, there exist a constant $\delta > 0$ and a sequence $(x_j) \subset \Omega$ such that $x_j \to x_0$ as $j \rightarrow + \infty$ and $B(x_j, \delta |x_j - x_0|^\alpha) \subseteq \Omega$.

Let $1 < p < n$. Consider a cutoff function $\eta \in C^\infty_0(B(0, \delta))$ with $\eta = 1$ around $0$ and extremal functions for $K(n,p)$ to one-parameter $\eps > 0$ given by
\[
v_\eps (x) = \eps^{-n/p^*} v_1\left(\frac{x}{\eps}\right)\, ,
\]
where $v_1$ is a extremal function satisfying $\| v_1 \|_{p^*} = 1$, see \cite{Au} or \cite{Ta}.

Define $w_j(x) = \eta(\frac {x-x_j}{\eps_j^\alpha}) v_{\eps_j^\beta}(x-x_j)$, where $\eps_j = |x_j-x_0|$ and $\beta > \alpha$ is a fixed number. We now compute all integrals involving the functions $w_j$ required in the proof of Theorems \ref{T.1} and \ref{T.2}.

Our estimates are simplified by using a trick of reduction to a previously known situation. Indeed, let $z_{\eps}(x) = \eta(x) v_\eps(x)$. Then, we can write

\[
w_j(x) = \eta\left(\frac {x-x_j}{\eps_j^\alpha}\right) v_{\eps_j^\beta}(x-x_j) = (\eps_j^{\alpha})^{-n/p^*} \eta\left(\frac {x-x_j}{\eps_j^\alpha}\right) (\eps_j^{\beta - \alpha})^{-n/p^*} v_1\left(\frac {x-x_j}{\eps_j^\alpha}/\eps_j^{\beta - \alpha}\right)
\]

\[
= (\eps_j^{\alpha})^{-n/p^*} z_{\eps_j^{\beta - \alpha}}\left(\frac {x-x_j}{\eps_j^\alpha}\right)\, ,
\]
and so,

\[
\nabla w_j(x) = \eps_j^{-n\alpha/p} \nabla z_{\eps_j^{\beta - \alpha}}\left(\frac {x-x_j}{\eps_j^\alpha}\right)\, .
\]
Therefore,

\[
\int_{\R^n} |\nabla w_j|^p\; dx = \int_{\R^n} |\nabla z_{\eps_j^{\beta - \alpha}}|^p\; dx,
\]

\[
\int_{\R^n} |w_j|^{p*}\; dx = \int_{\R^n} |z_{\eps_j^{\beta - \alpha}}|^{p*}\; dx
\]

\n and

\[
\int_{\R^n} |w_j|^{p}\; dx = \eps_j^{p \alpha} \int_{\R^n} |z_{\eps_j^{\beta - \alpha}}|^{p}\; dx
\]
Using then the corresponding estimates for $z_\eps$ obtained by Garcia Azorero and Peral Alonso in \cite{AzPe} in the right-hand side of the above inequalities, we get

\[
\int_{\R^n} |\nabla w_j|^p\; dx = K(n,p)^{-p} + O(\eps_j^{(n-p)(\beta - \alpha)/(p-1)}),
\]

\[
\int_{\R^n} |w_j|^{p*}\; dx = 1 + O(\eps_j^{n(\beta - \alpha)/(p-1)})
\]
and, for $n > p^2$,

\[
\int_{\R^n} |w_j|^{p}\; dx = \eps_j^{p \alpha} \left( a \eps_j^{p(\beta - \alpha)}+ O(\eps_j^{(n-p)(\beta - \alpha)/(p-1)}) \right)
\]

\[
= a \eps_j^{p\beta}+ O(\eps_j^{p \beta + (n-p^2)(\beta - \alpha)/(p-1)})\, ,
\]
where

\[
a = \int_{\R^n} |v_1|^p\; dx\, .
\]
Finally, given $\theta > 0$ and using that $\alpha \geq 1$, we derive

\[
\int_{\R^n} |x - x_0|^\theta |\nabla w_j|^p\; dx = O(\eps_j^{\theta})\, ,
\]
since the triangular inequality gives

\[
|x - x_0| \leq |x - x_j| + |x_j - x_0| \leq \delta \eps_j^\alpha + \eps_j = O(\eps_j)\, .
\]
Note that the last estimate is sharp because we have the reverse inequality

\[
|x - x_0| \geq |x_j - x_0| - |x - x_j| \geq \eps_j - \delta \eps_j^\alpha = O(\eps_j)\, .
\]
So, we end this section.

\section{Proof of Theorems \ref{T.1} and \ref{T.2}}

Assume that $A(x)$ and $a(x)$ satisfy the assumptions of Theorems \ref{T.1} and \ref{T.2}, respectively. According to the discussion done in the introduction on the proof of these results, it suffices to prove the strict inequalities (\ref{M1}) and (\ref{M2}).

By (\ref{H1}), there exists a constant $R > 0$ such that

\[
c_A \leq \inf_{u \in W^{1,2}_0(\Omega \cap B_R) \setminus \{0\}}
\frac{\int_{\Omega \cap B_R} \nabla^t u A(x) \nabla u\; dx -
\int_{\Omega \cap B_R} \lambda u^2\; dx}{\left(
\int_{\Omega \cap B_R}|u|^{2n/(n-2)} \;dx \right)^{(n-2)/n}}
\]

\[
\leq \inf_{u \in W^{1,2}_0(\Omega \cap B_R) \setminus \{0\}}
\frac{\int_{\Omega \cap B_R} \nabla^t u A(x_0) \nabla u\; dx + C_0
\int_{\Omega \cap B_R} |x - x_0|^{\gamma} |\nabla u|^2\; dx - \lambda
\int_{\Omega \cap B_R} u^2\; dx}{\left( \int_{\Omega \cap B_R}|u|^{2^*}
\;dx \right)^{2/2^*}}\, ,
\]
since $W^{1,2}_0(\Omega \cap B_R) \subset W^{1,2}_0(\Omega)$, where $B_R$ denotes the ball of radius $R$ centered at $x_0$.

In order to unify the proofs, we apply the change of variable $y = D P x$ in the right-hand side of the above inequality, where $P$ is an orthogonal matrix such that $P A(x_0) P^t$ is a diagonal matrix with entries $\lambda_i$ and $D$ is the diagonal matrix with entries $\lambda_i^{-1/2}$. For convenience, we denote $D P(\Omega \cap B_R)$ and $D P x_0$, respectively, by $\Omega \cap B_R$ and $x_0$. So, we can estimate $c_A$ as

\[
c_A \leq \inf_{u \in W^{1,2}_0(\Omega \cap B_R) \setminus \{0\}}
\frac{m_A^{1/n} \int_{\Omega \cap B_R} |\nabla u|^2\; dx + C_1 \int_{\Omega \cap B_R} |x - x_0|^{\gamma} |\nabla u|^2\; dx -
\lambda C_2 \int_{\Omega \cap B_R} u^2\; dx}{\left(
\int_{\Omega \cap B_R}|u|^{2n/(n-2)} \;dx \right)^{(n-2)/n}}\; ,
\]
where $C_1$ and $C_2$ are positive constants. Since affine maps preserve the order of singularity, the boundary of $\Omega \cap B_R$ is $\alpha$-singular at the new point $x_0$.

Similarly, by the assumption (\ref{H2}), we also derive

\begin{equation} \label{min}
c_{p,a} \leq \inf_{u \in W^{1,p}_0(\Omega \cap B_R) \setminus \{0\}}
\frac{m_{a,p} \int_{\Omega \cap B_R} |\nabla u|^p\; dx + C_0
\int_{\Omega \cap B_R} |x - x_0|^\sigma |\nabla u|^p\; dx - \lambda
\int_{\Omega \cap B_R} |u|^p\; dx}{\left( \int_{\Omega \cap B_R}|u|^{p^*}
\;dx \right)^{p/p^*}}\; .
\end{equation}
So, it suffices to prove only the strict inequality (\ref{M2}). By assumption, we have $\alpha \in [1, \sigma \frac{n-p^2}{np-p^2})$, so that we can choose a fixed number $\beta \in (\alpha \frac{n-p}{n-p^2}, \frac \sigma p)$. Thus, we get

\[
1 \leq \alpha < \beta,
\]

\[
p \beta < \sigma
\]
and

\[
p \beta < \frac{n-p}{p-1} (\beta - \alpha) < \frac{n}{p-1} (\beta - \alpha)\, .
\]

Finally, we compute the quotient of the right-hand side of (\ref{min}) for $u = w_j$ defined in the previous section. Note that $w_j \in W^{1,p}_0(\Omega \cap B_R)$ for $j$ large enough. So, using the integral estimates obtained for $w_j$ and the above inequalities for $\alpha$ and $\beta$, we deduce that

\[
c_{p,a} \leq \frac{m_{a,p} K(n,p)^{-p} - a \lambda \varepsilon_j^{p \beta} +
o(\varepsilon_j^{p \beta})}{1 + O(\eps_j^{n(\beta - \alpha)/(p-1)})} < m_a
K(n,p)^{-p}
\]

\n for each fixed $\lambda > 0$, provided $j$ is large enough. This ends the proof of both theorems.

\end{document}